\documentclass{article}
\usepackage{latexsym,amssymb,amsmath}
\usepackage{amsfonts}

\newtheorem{theorem}{Theorem}[section]

\newtheorem{proposition}[theorem]{Proposition}
\newtheorem{corollary}[theorem]{Corollary}
\newtheorem{definition}[theorem]{Definition}
\newtheorem{remark}[theorem]{Remark}
\newenvironment{proof}[1][Proof]{\textit{#1.} }{\ \rule{0.5em}{0.5em}}

\newcommand{\eps}{\varepsilon}

\numberwithin{equation}{section}

\begin{document}

\title{\bf Invariant measures for a stochastic Kuramoto-Sivashinky equation}
\author{{\bf Benedetta Ferrario}\footnote{B. Ferrario, Dipartimento di
      Matematica ''F. Casorati'', Universit\`a degli Studi di Pavia,
via Ferrata 1, 27100 Pavia, Italy; E-mail: benedetta.ferrario@unipv.it.}\\
         Dipartimento di Matematica, Universit\`a di Pavia\\
         Pavia, Italy}
        \date{ }
\maketitle

\noindent
{\bf Abstract:} 
For the 1-dimensional Kuramoto--Sivashinsky equation with random 
forcing term,
existence and uniqueness of 
solutions is proved.
Then, the Markovian semigroup is well defined; its properties are
analyzed, in order to provide sufficient conditions for 
existence and uniqueness of invariant
measures for this stochastic equation.
Finally, regularity results are presented.
\\[2mm]
\\[2mm]

\section{INTRODUCTION}\label{intro}
We consider the 1-dimensional Kuramoto--Sivashinsky equation perturbed
by an additive noise:
\begin{equation}\label{ks}
 du(t,x)+[\nu u_{xxxx}(t,x)+u_{xx}(t,x)+u(t,x)u_x(t,x)]\ dt=dW(t,x)
\end{equation}
where $t \in [0,T]$, $x \in \mathbb R$;
$\nu>0$ is a given coefficient. 
By $u_x,u_{xx},u_{xxxx}$ we denote, respectively,  the first, second 
and fourth derivative of $u$ with
respect to the space variable $x$.
Periodic conditions (with period $L$)
are assumed and an initial data $u(0,x)$ is
assigned.
In the right hand side of equation \eqref{ks} there is 
a Wiener process $W$ with covariance 
$
 \mathbb E [W(t,x)W(t^\prime,x^\prime)]=(t \wedge t^\prime) q(x,x^\prime)
$.

This stochastic equation is presented in the physical literature 
(see \cite{km,cl,uso} and references therein) in relation to
 a model for erosion by ion sputtering and 
 has been studied from the mathematical point of view in \cite{de} and
\cite{dy}.
\\
The deterministic Kuramoto--Sivashinsky equation, i.e. equation
\eqref{ks} without noise, has been 
introduced in \cite{k,kt,s1,sm1} in the 70's and from that time has
attracted the interest of many mathematicians
(see, e.g., \cite{t-inf}
for the basic results, and the references therein). 
It is known that it has a finite-dimensional
maximal attractor and an inertial manifold.
Numerical studies show chaotic behavior of its solutions, see \cite{k}.
However, it may happen that the dynamics has a more regular behavior
as far as statistical quantities, i.e. ensemble averages, 
are involved. 
This happens, for instance, 
in fluid dynamics: the individual solutions are chaotic but
statistical properties of the dynamics are more regular, as
investigated by  turbulence theory.
The results on invariant measures, presented in this paper for the 
stochastic Kuramoto--Sivashinsky equation \eqref{ks},
are exactly about the statistical behavior of the solution.

Comparing our results with those of \cite{de} and \cite{dy}, we notice that
we construct  solutions to equation \eqref{ks} under 
assumptions on the noise weaker
than in \cite{de} and \cite{dy}
so to analyze the case presented in \cite{km} and \cite{cl}; moreover
we deal with invariant measures for equation \eqref{ks}.
Existence of invariant measures may be obtained from the results in 
\cite{dy}; indeed, the estimates to prove existence of 
a finite dimensional random attractor are similar to the ones used in
 \cite{fla} for existence of an invariant measure.
But we prove it by another technique and with
weaker assumptions.
Moreover we do not restrict to the case of odd solutions.
Finally, we tackle the problem of uniqueness of invariant measures.

We notice that the deterministic equation
is a fourth-order PDE with the non linear term of 
the same form as in the Burgers or one-dimensional Navier-Stokes equations.
We shall exploit this peculiarity in this paper, borrowing some
techniques
used in the analysis of the stochastic Burgers  or 
Navier--Stokes equations.

This  article is organized as follows.
In Section \ref{intro} we introduce the stochastic
Kuramoto--Sivashinsky equation as an It\^o equation in Hilbert spaces.
In Section \ref{soluts} we prove a theorem of existence and uniqueness
of solutions.
In Section \ref{invmeas} we explain how to prove existence of invariant
measures (the details are given in Section \ref{exinv}) and uniqueness
(the details are given in Sections \ref{irred} and \ref{strongly}).
Regularity results  are proved in Section
\ref{secGIR}; in this way the results of irreducibility and strongly
Feller property, proved first  in the basic space $H$,
 are extended to more regular spaces.
Section \ref{conclu} 
presents the final theorem, covering all the results proved.

\section{ABSTRACT SETTING}
In this section, we  introduce 
spaces and operators in order to define 
an abstract formulation of equation \eqref{ks}
as an It\^o equation in Hilbert spaces.

\medskip
\centerline{SPACES }

Let $\mathcal P$ be the space of periodic  $C^\infty$-functions
defined on $[-\frac L2,\frac L2]$ and with zero mean.
Closing this space with respect to the 
$L^2(-\tfrac L2, \tfrac L2)$ and $H^2(-\tfrac L2, \tfrac L2)$-norm 
we get the following spaces:
$$
H=\{u \in L^2(-\tfrac L2, \tfrac L2): \int_{-L/2}^{L/2} u(x)\ dx=0\},
$$
a Hilbert space with scalar product $\langle u, v \rangle_H=
\int_{-L/2}^{L/2}u(x) v(x) \ dx$,
\\
and
$$
V=H\cap \{u \in H^2(-\tfrac L2, \tfrac L2):u(-\tfrac L2)=u(\tfrac L2),
u_x(-\tfrac L2)=u_x(\tfrac L2)\},
$$
a Hilbert space with scalar product $\langle u, v \rangle_V=
\int_{-L/2}^{L/2}u_{xx}(x) v_{xx}(x) \ dx$.

\medskip

\centerline{OPERATORS}

We define the operator $A$ as
$$
 Au=-u_{xx}, \quad D(A)=V.
$$
This is a linear operator in $H$, densely defined.
It is strictly  positive; its eigenvectors and eigenvalues
are
$$
 \tilde e_{j,1}(x)=\sqrt{\tfrac 2L} \sin(\tfrac{2j\pi}L x), \;
 \tilde e_{j,2}(x)= \sqrt{\tfrac 2L} \cos(\tfrac{2j\pi}L x),
$$
$$
 \tilde \lambda_j=\frac{4 \pi^2}{L^2} j^2
$$
for $j=1,2, \dots$.
To shorten notations, from now on we shall denote by 
$\{e_j\}_{j=1}^\infty $ the sequence of the eigenvalues with
corresponding eigenvectors $\lambda_j$ (this is nothing but a 
relabelling of the sequence: $e_{2k-1}=\tilde e_{k,1}$, 
$e_{2k}=\tilde e_{k,2}$ 
and $\lambda_{2k-1}=\lambda_{2k}=\tilde \lambda_k$ for $k=1,2,\dots$).
The sequence of the eigenvectors of $A$ is a complete orthonormal basis
of the space $H$. This implies that every $u\in H$ can be written as
$u=\sum_j u_j e_j$,
where the coefficients $u_j$ satisfy the condition
$\sum_j |u_j|^2 <\infty$.

The power operator
$A^\alpha$ exists for any $\alpha \in \mathbb R$
(see, e.g., \cite{pazy}), $D(A^\alpha)=\{u=\sum_j u_je_j: 
\sum_j \lambda_j^{2\alpha}  u_j^2<\infty\}$,
$A^\alpha u=\sum_j \lambda_j^\alpha
u_j e_j$ and $|A^\alpha u|^2=\sum_j \lambda_j^{2\alpha}  u_j^2$. 
Given any $\alpha$ and  $\beta$, the operator $A^\alpha$ is an isomorphism
from $D(A^\beta)$ to $D(A^{\beta-\alpha})$.
\\
We have that $V=D(A)$ and the Poincar\'e inequality
$$
 |u|_V\ge \lambda_1 |u|_H.
$$
Moreover, $D(A^2)=V\cap \{u \in H^4(-\frac L2, \frac L2):
u_{xx}(-\tfrac L2)=u_{xx}(\tfrac L2),
u_{xxx}(-\tfrac L2)=u_{xxx}(\tfrac L2)
\}$
and
$(D(A))^\prime=D(A^{-1})$,
where $(D(A))^\prime$ 
is the dual space of $D(A)$ with respect to the duality of the
$H$-scalar product.
For any natural $m\in \mathbb N$, the $D(A^{\frac m2})$-norm is
equivalent
to the $H^m(-\tfrac L2, \tfrac L2)$-norm.

In particular, for $\alpha>0$, $A^\alpha$ is a self-adjoint
operator in $H$ generating a $C_0$-semigroup of linear operators in $H$: 
$$
 e^{-A^\alpha t}u=\sum_je^{-\lambda_j^\alpha t}u_j e_j \;\;
\mbox{ for any }t\ge 0, u =\sum_j u_j e_j \in H .
$$
In the following we shall use this result
(see, e.g., \cite{pazy}):
for any $\beta >0$
there exists a constant $M_\beta$ such that
\begin{equation}\label{regeA}
 |A^{2\beta} e^{-A^2 t}u|_H\le \frac {M_\beta} {t^\beta} |u|_H 
   \quad  \mbox{ for any } t>0, u\in H .
\end{equation}

The bilinear operator $B$
is studied investigating the associated 
trilinear form 
$$
 b(u,v,z)=\langle B(u,v),z\rangle_H 
          = \int_{-L/2}^{L/2} u(x) v_x(x) z(x)\ dx .
$$
Using H\"older inequality and the continuous 
embedding of spaces $H^1(-\frac L2, \frac L2)\subset
L^\infty(-\frac L2, \frac L2)$, we obtain that
there exists a constant $c$ such that
\begin{equation}\label{buzu}
|b(u_1,u_2,u_3)|\le |u_1|_H |(u_2)_x|_{L^\infty}|u_3|_H 
\le c |u_1|_H |Au_2|_H |u_3|_H 
\end{equation}
for  $u_1,u_2,u_3 \in \mathcal P$.
By density, the estimates hold for all $u_1,u_3 \in H$ and $u_2 \in D(A)$.

We shall use the following identities for elements of $\mathcal P$, 
obtained integrating by parts:
\begin{equation}\label{scambi}
\left[
\begin{array}{l}
 b(u,u,u)=0,
 \\
 b(u_1,u_2,u_2)=b(u_2,u_2,u_1)=-\frac 12 b(u_2,u_1,u_2),
\\
 b(u_1,u_2,u_3)=-b(u_2,u_1,u_3)-b(u_1,u_3,u_2).
\end{array}\right.
\end{equation}

We collect here useful estimates. Hereafter 
the symbol $c$ denotes different constants.
\begin{proposition}
\begin{align}
&|B(u,v)|_{D(A^{-1})}\le c |Au|_H |v|_H \label{uno}
\\\label{due}
&|B(z,v)|_{D(A^{-1})}\le  c |z|_H |Av|_H
\\\label{tre}
&|B(z,z)|_{D(A^{-1})}\le c |z|^2_H 
\\\label{sei}
&|B(u,v)|_{D(A^{-\delta})}
\le c |u|_{D(A^{\frac 12 -\delta })} |v|_{D(A^{\frac 12 -\delta })}
  \quad \text{ if } \;\delta \le 0 
\end{align}
\end{proposition}
\proof 
First we show the inequality for elements of $\mathcal P$; then 
by density they hold true for elements in the spaces specified by the
norms involved at each instance. 
We shall use  repeatedly  H\"older inequality and the fact
that $H^1(-\frac L2,\frac L2)$ 
is continuously embedded in $L^\infty(-\frac L2,\frac L2)$, so that
$|u|_{L^\infty}\le c|u_x|_H$. Moreover, $|B(u,v)|_{D(A^{-\alpha})}
=\displaystyle\sup_{| w|_{D(A^\alpha)}\le 1}|b(u,v,w)|$.
\\
For each inequality we show the main estimates  to get it.

\eqref{uno}: we have $|B(u,v)|_{D(A^{-1})}\le c |A^{1/2}u|_H |v|_H$, since
$|b(u,v,w)|\le |b(v,u,w)|+|b(u,w,v)|$ by \eqref{scambi}
 and
 $|b(v,u,w)|\le
|v|_H |u_x|_H |w|_{L^\infty}$, $|b(u,w,v)|\le |u|_{L^\infty} |w_x|_H |v|_H$.

\eqref{due}: $|B(z,v)|_{D(A^{-1})}\le \frac 1{\lambda_1}|B(z,v)|_H 
     \le \frac 1{\lambda_1} |z|_H |v_x|_{L^\infty}$.

\eqref{tre}: $|b(z,z,w)|=\frac 12 |b(z,w,z)|$ by \eqref{scambi} and 
             $|b(z,w,z)|\le   |z|^2_H |w_x|_{L^\infty}$.


\eqref{sei}: we have $|u v_x|_H \le |u|_{L^\infty} |v_x|_H
\le c |u|_{D(A^{1/2})} |v|_{D(A^{1/2})}$ and
$|(u v_x)_x|_{H}
\le |u_x v_x|_H+|u v_{xx}|_H
\le |u_x|_H |v_x|_{L^\infty}+|u|_{L^\infty} |v_{xx}|_H
\le c |u|_{D(A^{1/2})} |v|_{D(A)}$. Hence
\begin{align*}
|B(u,v)|_H&\le c|u|_{D(A^{1/2})} |v|_{D(A^{1/2})}
\intertext{and}
|B(u,v)|_{D(A^{1/2})} &\le c |u|_{D(A^{1/2})} |v|_{D(A)}.
\end{align*}
By interpolation 
$$|B(u,v)|_{D(A^\theta)}\le c
|u|_{D(A^{1/2})} |v|_{D(A^{\frac 12 +\theta })}\, , \qquad 0<\theta<\frac 12.
$$
We conclude the case of $-\frac 12 \le \delta <0$ 
noting that $|u|_{D(A^{\frac 12})}\le c_\delta |u|_{D(A^{\frac 12-\delta})}$.
\\
For $\delta <-\frac 12$, set $\frac m2=-\delta$.
$D(A^{\frac m2})$ is a multiplicative algebra for $m=1,2,\dots$.
Then the result is trivial for $m$ integer;
the estimate is even better:
$$
 |B(u,v)|_{D(A^{\frac m2})} = |u v_x|_{D(A^{\frac m2})}
 \le c |u|_{D(A^{\frac m2})} |v|_{D(A^{\frac 12 +\frac m2})} .
$$
This allows to extend the result to any $m > 1$ as 
before by interpolation. \hfill
$\Box$

\smallskip
To shorten notations, we write $B(u)$ for $B(u,u)$.

\medskip
\centerline{ABSTRACT FORMULATION}

The abstract formulation of the initial value problem 
for  equation \eqref{ks} is
\begin{equation}\label{sks}
\left\{
\begin{array}{l}
du(t)+ [ \nu A^2 u(t)-Au(t)+B(u(t))]dt=Gdw(t) , \\
u(0)=y.
\end{array}
\right.
\end{equation}

We have written the Wiener process as
$Gw(t)$, where $G$ is a linear operator  
and  $w$ is a cylindrical Wiener process in $H$
defined on a probability space with filtration
$(\Omega,\mathcal F,\{\mathcal
F_t\}_{t\ge 0},\mathbb P)$
(i.e. given a sequence $\{\beta_j\}_{j=1}$ of i.i.d. 
standard Wiener processes, we represent the Wiener process in series
as $w(t)= \sum_j \beta_j(t) e_j$).

Hereafter, 
we denote by $|\cdot|$ the $H$-norm and by
$\langle\cdot,\cdot\rangle$ the scalar product in $H$.
If other norms are involved,  they will be specified at each
instance.
 \\
$\mathbb E$ denotes the expectation with respect to the probability measure
$\mathbb P$.

\section{SOLUTION: EXISTENCE, UNIQUENESS AND PROPERTIES}
\label{soluts}
From now on, we assume that $G$ is a linear operator in
$\mathcal P$ such that
\begin{equation}\label{ipotG}
 \sum_{j=1}^\infty \lambda_j^{-2}\sum_{k=1}^\infty 
 [\langle Ge_k,e_j\rangle]^2<\infty.
\end{equation}
When possible, the operator $G$ in defined as a linear operator in $H$ 
(notice that $\mathcal P$ is dense in $H$) and this extension
is again denoted by $G$.
Then, condition \eqref{ipotG} is equivalent to require $A^{-1}G$  be 
a Hilbert--Schmidt operator in $H$. We point out that \cite{de} and 
\cite{dy} assume $G$ be a Hilbert--Schmidt operator, which is stronger
than \eqref{ipotG}.

\begin{remark}\label{convergenza}
Keeping in mind the expression of the eigenvalues, we find that
condition \eqref{ipotG} 
 is satisfied for instance if $G=LA^\gamma$ with $\gamma<\frac 34$
and $L$ an isomorphism in $H$ (we mean that $L$ acts on the basis
$\{e_j\}$ renaming it).
Actually, for $0<\gamma<\frac 34$, $G=LA^\gamma$  is a  linear operator in $H$
with domain $D(A^\gamma)$, whereas, for $\gamma\le 0$, $G=LA^\gamma$ 
is a  linear bounded operator in $H$.
Usually the literature on stochastic equations in Hilbert spaces
treats the case of $G$ linear bounded operator
in $H$; in this paper we assume $A^{-1}G$ be a linear bounded operator
in $H$ of Hilbert--Schmidt type.
\end{remark}

We point out that, in the physical literature, 
the noisy Kuramoto-Sivashinsky equation 
is
presented as
$$
  \partial_t h(t,x)+\nu h_{xxxx}(t,x)+h_{xx}(t,x)-\tfrac 12  h_x(t,x)^2\ 
  =\eta(t,x) ,
$$
where $\nu$ is a positive surface diffusion coefficient and the
variable $h$ is the height profile of a surface eroded by ion sputtering.
Setting $u=-h_x$ and $\partial_t w=-\eta_x$, we get \eqref{ks} and 
the unknown $u$ 
can be interpreted as a one-dimensional velocity field
in a compressible fluid (see \cite{uso}).

\cite{cl} and
\cite{km} consider a centered Gaussian noise $\eta$ with covariance 
$$
\mathbb E[\eta(t,x)\eta(t^\prime,x^\prime)]
=\delta(t-t^\prime) \delta(x-x^\prime),
$$
i.e. $\eta(t,x)=\sum_j \dot \beta_j(t) e_j(x)$.
This corresponds to have
$$
 Gw=\sum_{j \text{ even}} \beta_j \lambda_j^{1/2} e_{j-1}
-
     \sum_{j \text{ odd}} \beta_j \lambda_j^{1/2} e_{j+1} ,
$$
in our abstract equation \eqref{sks}, i.e.
$G=L A^{1/2}$  where $L$ is the  linear
bounded operator in $H$ defined by
$Le_j=(-1)^{j} e_{j+(-1)^{j+1}}$.
Actually, $L$ is an isomorphism in $H$. 
Then, \eqref{ipotG} is satisfied.

On the other hand, 
\cite{uso} considers a noise with covariance 
$$\mathbb E[\eta(t,x)\eta(t^\prime,x^\prime)]
=\delta(t-t^\prime) [I+A]\delta(x-x^\prime),
$$
i.e. $\eta(t,x)=\sum_j \dot \beta_j(t) [I+A]^{1/2}e_j(x)$.
This corresponds to have
$G=L A^{1/2}[I+A]^{1/2}$ in our equation \eqref{sks}. In this case
\eqref{ipotG} is not satisfied. 

\medskip

We begin defining what is a solution for equation \eqref{sks}
in this work. Given an initial data $y \in H$, 
we consider solutions with less spatial regularity than
whose of \cite{de}, \cite{dy}.
\begin{definition}\label{defin}
A stochastic process $u$ is 
a weak solution of equation \eqref{sks} on the time interval $[0,T]$
if, $\mathbb P$-a.s.
$$
 u \in C([0,T];H),
$$
it is progressively measurable and
it satisfies the following identity
\begin{equation}\label{def-u}
 \langle u(t),h\rangle+\int_0^t \langle u(s),\nu A^2 h
 - Ah\rangle ds -\frac 12 \int_0^t b(u(s),h,u(s))
 ds =
 \langle y,h\rangle + \langle h, Gw(t)\rangle
\end{equation}
for any $t \in [0,T]$ and $h \in \mathcal P$.
\end{definition}

Relationship \eqref{def-u} is formally
obtained from \eqref{sks} multiplying by a test function $h$
and
integrating over the spatial domain; integration by parts 
as in \eqref{scambi} yields 
the above expression. 
Moreover, elementary calculus based on the It\^o isometry gives
$$
 \mathbb E [\langle h,Gw(t)\rangle]^2 \le \;t \;
  \big(\sum_k \lambda_k^2 h_k^2\big)
 \sum_j \lambda_j^{-2}\sum_{k=1}^\infty 
 [\langle Ge_k,e_j\rangle]^2.
$$
Bearing in mind \eqref{ipotG} and \eqref{buzu},
we have that 
all the terms in \eqref{def-u} are well defined.

Analyzing the equation for $u$, we first study the linear part.
Notice that
the linear operator $\nu A^2-A $ is not strictly positive
(this depends on $L$ and $\nu$, because its eigenvalues are 
$\nu \lambda_j^2-\lambda_j$)
and the negative eigenvalues cause instability for the linear
Kuramoto--Sivashinsky equation. 

Hence, we introduce the linear stochastic equation 
\begin{equation}\label{eqOU}
 dz_a(t)+\nu A^2 z_a(t)\ dt -A z_a(t)\ dt
      +a z_a(t)\ dt = Gdw(t), \quad z_a(0)=\zeta.
\end{equation}
We fixe a value  $a>\frac 1{4\nu}$ so to have $\nu \lambda_j^2-\lambda_j+a>0$ 
for all $j$. Therefore the operator $\nu A^2-A+aI$, with domain $D(A^2)$, 
is strictly positive.
\\
The process solving \eqref{eqOU} is
$$
z_a(t)=e^{-(\nu A^2-A+a)t} \zeta+
\int_0^t e^{-(\nu A^2-A+a)(t-s)}Gdw(s) .
$$
Writing $e^{-(\nu A^2-A+a) (t-s) }Gw(s)
=\displaystyle\sum_{j,k}e^{-(\nu  \lambda_j^2- \lambda_j+a) (t-s)}
\langle Ge_k,e_j\rangle \beta_k(s) e_j$,
we have that 
\begin{equation}\label{media-z}
\begin{split}
\mathbb E[z_a(t)]&=e^{-(\nu A^2-A+a)t} \zeta \in H ,
\\
\\
 \mathbb E \left[|z_a(t)-\mathbb E z_a(t)|^2\right]
&=\sum_{j,k}  \int_0^t e^{-2(\nu\lambda_j^2- \lambda_j+a)(t-s)}|
          \langle G e_k,e_j\rangle |^2 ds
\\
& \le 
  \sum_{j,k} |\langle G e_k,e_j\rangle |^2
       \frac 1{2(\nu \lambda_j^2- \lambda_j+a)} .
\end{split}
\end{equation}
Therefore, if condition \eqref{ipotG} holds, then 
for any time $t$, the Gaussian  random variable $z_a(t) \in H$, $
\mathbb P$-a.s.
Using the factorization method as in \cite{dpz}, we get also that
there exists a continuous version of $z_a$ with values in $H$ (and 
from now on we
shall consider  this continuous version).

\bigskip

To solve  equation \eqref{sks}, we introduce a new unknown $v_a$
as suggested in analysis of other  It\^o equations with additive noise
(see, e.g., \cite{fla}). 
Set $v_a=u-z_a$ and $z_a(0)=0$. 
Making the difference of the equations satisfied by $u$ and $z_a$,
we obtain that  
the equation satisfied by $v_a$ does not contain the 
noise term. This is
\begin{equation} \label{eq-v}
\left\{
\begin{array}{l}
 \frac{d}{dt} v_a(t) +\nu A^2 v_a(t) - A v_a(t) +B(v_a(t)+z_a(t))=a z_a(t),
 \\
 v_a(0)=y.
 \end{array}\right.
\end{equation}
We first prove a result for this problem, considered pathwise. This
resembles the result for the deterministic equation (see \cite{t-inf}).
\begin{proposition}\label{vinH}
Assume \eqref{ipotG}.
Then, for any $y \in H$  and $T>0$ there exists a unique solution $v_a$ for
\eqref{eq-v}
such that
$$
 v_a \in C([0,T];H) \cap L^2(0,T;D(A))
$$
$\mathbb P$-a.s.
\end{proposition}
\proof
We prove existence of a solution by means of the Galerkin method,
i.e. we first deal with a finite dimensional problem for which there
exists a solution and then we pass to the limit to recover the
original evolutionary problem.

For any $n \in \mathbb N$, let $\Pi_n$ be the orthogonal projector from $H$ to 
the space spanned by $e_1,e_2,\dots,e_n$.
Set $v_n=\Pi_n v$, $B_n=\Pi_n B$. Notice that the operators $A$ and
$\Pi_n$ commute.
The Galerkin system is
\begin{equation}\label{Gal-v}
 \left\{
\begin{array}{l}
 \frac{d}{dt} v_a^n(t) + \nu A^2 v_a^n(t) - A v_a^n(t)  + 
B_n(v_a^n(t)+z_a^n(t))=  a z_a^n(t),
 \\
 v_a^n(0)=\Pi_n y.
 \end{array}\right.
\end{equation}
We work pathwise.
Since the coefficients are locally Lipschitz, 
there exists a unique solution, local in time. To show global
existence we need a priori estimates.
We take the scalar product of this equation with $v_a^n$,
use \eqref{scambi}, \eqref{buzu} and Young inequality, to  obtain
\[
\begin{split}
 \frac 12 \frac d{dt} |v_a^n|^2 +\nu |Av_a^n|^2 
&= \langle Av_a^n,v_a^n\rangle
  +a \langle z_a^n, v_a^n\rangle 
   -b(v_a^n,v_a^n,v_a^n)-b(v_a^n,z_a^n,v_a^n)
\\
&\qquad -b(z_a^n,v_a^n,v_a^n)-b(z_a^n,z_a^n,v_a^n)
\\
&=
 \langle Av_a^n,v_a^n\rangle
+ a \langle z_a^n, v_a^n\rangle +b(z_a^n,v_a^n,v_a^n)
        +\tfrac 12 b(z_a^n,v_a^n,z_a^n)
\\
& \le
|v_a^n| |Av_a^n| +a |z_a^n| |v_a^n| + c |z_a^n| |v_a^n| |Av_a^n| 
  +\tfrac 12 c |z_a^n|^2 |Av_a^n| 
\\
& \le
\tfrac 12 \nu|Av_a^n|^2 + c_\nu (1+|z_a^n|^2)|v_a^n|^2 
     +\frac {a^2}2 |z_a^n|^2+ c_\nu |z_a^n|^4.
\end{split}
\]
$c_\nu$ denotes different constants depending on $\nu$ and on the
spatial domain.
Then, there exists a constant $C_1$ (depending on $\nu$ and on the
lenght $L$, but not on $n$) such that 
\begin{equation} \label{stime-v}
   \frac d{dt} |v_a^n|^2 +\nu |Av_a^n|^2 \le
   C_1 (1+|z_a|^2)|v_a^n|^2 +  a^2 |z_a|^2+C_1 |z_a|^4 .
\end{equation}
Applying Gronwall inequality to $ \frac d{dt} |v_a^n|^2 \le
   C_1 (1+|z_a|^2)|v_a^n|^2 + a^2 |z_a|^2+C_1 |z_a|^4$
we get that
\begin{equation}\label{sup-v}
|v_a^n(t)|^2 
\le |y|^2 e^{\int_0^t C_1 (1+|z_a(s)|^2) ds} 
     +\int_0^t e^{\int_s^t C_1(1+|z_a(\tau)|^2) d\tau}
     [a^2|z_a(s)|^2+C_1|z_a(s)|^4]ds .
\end{equation}
Since $\displaystyle\sup_{0\le s \le T}|z_a(s)|$ is finite $\mathbb P$-a.s.
under assumption \eqref{ipotG},
we conclude that
$$
 \sup_{0 \le t \le T} |v_a^n(t)|^2\le C_2,
$$
where $C_2$ is independent of $n$. 
Now, integrating in time \eqref{stime-v}, by means of the last
estimate we get also that
$$
 \int_0^T |Av_a^n(t)|^2 dt \le C_3,
$$
where $C_3$ is independent of $n$. 
Hence, $v_a^n \in C([0,T];H)\cap L^2(0,T;D(A))$.

Finally, we note that
$$
 \frac{d}{dt} v_a^n =- \nu A^2 v_a^n 
+ A v_a^n + a z_a^n - B_n(v_a^n)-B_n(z_a^n)-B_n(v_a^n,z_a^n)-B_n(z_a^n,v_a^n).
$$
The r.h.s. belongs to $L^2(0,T;D(A^{-1}))$; this is easy to check for
the first three terms, whereas for the terms involving $B_n$ we have to
bring to mind  \eqref{uno}-\eqref{tre}.

Summing up, 
the Galerkin sequence $v_a^n$ is bounded in 
$L^\infty(0,T;H)\cap L^2(0,T;D(A))\cap H^1(0,T;D(A^{-1}))$.
Since the space $L^2(0,T;D(A))\cap H^1(0,T;D(A^{-1}))$ 
is compactly embedded in $L^2(0,T;H)$
(see, e.g., \cite{temam} at pg. 271), 
we conclude that there exists a subsequence
$v_a^{n_k}$  and a limit $v_a$ such that

$v_a^{n_k}$ converges to $v_a$ \; weakly in $L^2(0,T;D(A))$,
 
$v_a^{n_k}$ converges to $v_a$ \; $\star$-weakly in $L^\infty(0,T;H)$,

$v_a^{n_k}$ converges to $v_a$ \; strongly in $L^2(0,T;H)$.
\\
These convergences grant that $v_a$ is a solution
of equation \eqref{eq-v}; notice that the strong convergence allows to
pass to the limit in the non linear term $B$ (see details in
\cite{temam},
dealing with the Navier--Stokes equation which has this same non linearity).

Finally, if $v_a \in L^2(0,T;D(A))$, 
$v_a^\prime \in L^2(0,T;D(A^{-1}))$, 
then $v_a \in C([0,T];H)$ (see \cite{temam}, Chapter III Lemma 1.2).

Uniqueness is easy to check. We shall prove it for $u$ in the next
theorem and the method applies also to $v_a$ (usually it is more
difficult to get uniqueness for $u$, because $u$ is less regular than
$v_a$; for this reason we give
details only for $u$). \hfill
$\Box$

\medskip

The result for the process $u$ is the following.
\begin{theorem}\label{primoTEO}
Assume \eqref{ipotG}.
Then, for any $y \in H$ and $T>0$ there exists a unique solution $u $ for
\eqref{sks}
as defined in Definition \ref{defin} such that
$$
u \in C([0,T];H)
$$
$\mathbb P$-a.s..
Moreover, $u$ is a Markov process in $H$, which is Feller in $H$.
\end{theorem}
\proof
The process $u=v_a+z_a$ is a solution of \eqref{sks} by construction and
the regularity of $v_a$ and $z_a$ provides $u\in C([0,T];H) $, $\mathbb P$-a.s.

As far as we are concerned with the the uniqueness, 
let $u_1, u_2 \in C([0,T];H)$ be two solutions of equation \eqref{sks}
with the same initial data. Set $U=u_1-u_2$.
Then this difference satisfies the following equation
$$
 \frac {dU}{dt} +\nu A^2 U -AU +B(u_1)-B(u_2)=0;\qquad U(0)=0.
$$
By the bilinearity of $B$: $B(u_1)-B(u_2)=
B(U,u_1)+B(u_2,U)$. Then
$$
 \frac {dU}{dt} +\nu A^2 U = AU - B(U,u_1) - B(u_2,U) .
$$
Taking the scalar product in $H$ of  this equation  with $U$ and
proceeding to estimate the terms as done before for $v_a$, we get
\[
\begin{split}
 \frac 12 \frac d{dt}|U|^2+\nu |AU|^2 
&\le |U||AU|+c|U||AU|(|u_1|+|u_2|)
\\
&\le
    \frac \nu 2 |AU|^2+c_\nu (1+|u_1|^2+|u_2|^2)|U|^2.
\end{split}
\]
Actually, we should work first with the finite dimensional
approximation and then pass to the limit, obtaining that
$U \in C([0,T];H)\cap L^2(0,T;D(A))$. The fact that $U$ is more regular
than $u_1$ and $ u_2$ justifies all the estimates.

From \eqref{stimeU} we conclude by Gronwall lemma that
\begin{equation}\label{stimeU}
\sup_{t \in [0,T]}|U(t)|^2 \le |U(0)|^2 
e^{2c_\nu \int_0^T(1+|u_1(s)|^2+|u_2(s)|^2)ds} .
\end{equation}
Then
$$
 |U(t)|=0 \mbox{ for all } t \in [0,T]
$$
and uniqueness is proved.

Since the processes $v_a^n$ and $z_a^n$ are progressively measurable, so
are also $u^n$ and in the limit $u$ too.
The same for the Markov property, 
i.e. given $0<t_1<\dots<t_m\le T$
$$
 \mathbb P\{u(t_m)\in\Gamma|u(t_1),\dots,u(t_{m-1})\}
=
 \mathbb P\{u(t_m)\in\Gamma|u(t_{m-1})\}
$$
for any Borelian subset $\Gamma$ of $H$.

Finally, denoting by
$u(\cdot;y)$
the solution to equation \eqref{sks} with initial data $y\in H$, 
the Markovian transition semigroup $\{P_t\}_{t \ge 0}$, defined by 
$$
 (P_t\phi)(y)=\mathbb E \phi(u(t;y)) ,
$$
is well defined  on the space of Borelian bounded
functions $B_b(H)$.
Moreover, it is Feller, that is $P_t:C_b(H)\to C_b(H)$.
Indeed, from \eqref{stimeU}, we get that the solution $u$ depends continuously 
on the initial data; then for any $t$
$$
 \lim_{y_1 \to y_2}|u(t;y_1) -u(t;y_2)|=0 
 \quad \mathbb P -a.s. .
$$
Given any  continuous function $\phi:H\to \mathbb R$, we have that
 $\phi(u(t;y_1)) \to \phi(u(t;y_2))$ pathwise as $y_1 \to y_2$. 
Since $\phi$ is bounded, 
by dominated convergence theorem we conclude  that
$\mathbb E\phi(u(t;y_1)) \to \mathbb E\phi(u(t;y_2))$ as $y_1 \to y_2$.
\hfill
$\Box$

\medskip

If the process $z_a$, solution to the linear stochastic equation
 \eqref{eqOU},
 is more regular, then also
$u=v_a+z_a$ is more regular. We note that, by interpolation,  
$$
  |A^\alpha v|\le c |v|^{1-\alpha} |Av|^\alpha
$$
for $0<\alpha<1$.
Hence, if $v_a \in C([0,T];H)\cap L^2(0,T;D(A))$, we have
$v_a \in L^{2/\alpha}(0,T;D(A^\alpha))$.

Moreover if 
\begin{equation}\label{ipoalfa}
 \sum_{j=1}^\infty \lambda_j^{2(\alpha-1)}\sum_{k=1}^\infty 
 [\langle Ge_k,e_j\rangle]^2<\infty
\end{equation}
then $z_a \in C([0,T];D(A^\alpha))$ $\mathbb P$-a.s. .
Indeed the estimates to prove this  
are easily obtained as in \eqref{media-z}. 
We write them for a general initial data $\zeta \in D(A^\alpha)$:
\begin{equation}\label{zinalfa}
\begin{split}
 \mathbb E [z_a(t)]&=e^{-(\nu A^2-A+a)t} \zeta \in D(A^\alpha) ,
\\
\mathbb E \left[|z_a(t)-\mathbb E z_a(t)|_{D(A^\alpha)}^2\right]
& \le 
  \sum_{j,k} |\langle G e_k,e_j\rangle |^2
  \frac {\lambda_j^{2\alpha}}{2(\nu \lambda_j^2-\lambda_j+a)} .
\end{split}
\end{equation}
This implies the following result, which enforces 
the previous theorem.
\begin{corollary}\label{cor-alfa}
(i) If \eqref{ipoalfa} holds for some $0<\alpha <1$, then 
there exists a unique process $u$
 solution to \eqref{sks}, which in addition to the properties of 
 Theorem \ref{primoTEO} has 
$u \in C([0,T];H)\cap L^{2/\alpha}(0,T;D(A^\alpha))$ $\mathbb P$-a.s.

(ii)
If \eqref{ipoalfa} holds  for some $\alpha \ge 1$, then 
there exists a unique process $u$
solution to \eqref{sks}, which in addition to the properties of 
Theorem \ref{primoTEO} has
$u \in C([0,T];H)\cap L^2(0,T;D(A))$ $\mathbb P$-a.s.
\end{corollary}

\section{INVARIANT MEASURES}\label{invmeas}
\begin{definition}
A probability measure $\mu$  on $H$ is invariant for the Markovian semigroup 
$\{P_t\}_{t \ge 0}$ associated to equation \eqref{sks} if
$$
 \int P_t \phi d\mu = \int \phi d\mu \qquad \mbox{ for all } \phi\in
 C_b(H),\; t\ge 0.
$$
\end{definition}
Introducing the semigroup $\{P^*_t\}_{t\ge 0}$, acting on probability
measures 
on $H$, as
$\int \phi \, d P^*_t\mu=\int P_t\phi \, d\mu \equiv <P_t\phi,\mu>$,
a measure $\mu$ is invariant if
$$
   P^*_t \mu = \mu \mbox{ for all } t\ge 0,
$$
that is, $\mu $ is a fixed point for the evolution of probability measures
under $P_t^*$.

Since the Markovian semigroup $\{P_t\}_{t\ge 0}$ is Feller in $H$, we use the
well-known
Krylov--Bogoliubov method to prove existence of invariant measures.
Namely, if the family of measures
$$
\mu_{T}=\frac 1T \int_1^T P_s^* \delta_0 ds , \qquad T> 1 ,
$$
(where $P_s^* \delta_0$ is the law of the process $u(s;0)$)
is tight in $H$, then there exists a subsequence $\mu_{T_k}$
weakly convergent to a measure $\mu$, as $k \to\infty$
(and $T_k \to \infty$). Then this limit measure $\mu$ is
invariant for the Markovian semigroup $\{P_t\}_{t\ge 0}$; 
indeed, for any $\phi \in C_b(H)$
$$
\begin{array}{rl}
 <P_t\phi,\mu>
&={\displaystyle
 \lim_{k\to \infty}} \frac 1{T_k} \int_1^{T_k} <P_t \phi,P_s^* \delta_0> ds
\\
&={\displaystyle\lim_{k\to\infty}} 
       \frac 1{T_k} \int_1^{T_k} <\phi,P_{t+s}^* \delta_0> ds
\\
&={\displaystyle\lim_{k\to\infty}}
   \frac 1{T_k} \int_{t+1}^{t+T_k} <\phi,P_s^* \delta_0> ds 
\\
&={\displaystyle\lim_{k\to\infty}}
        \frac 1{T_k} \int_1^{T_k} <\phi,P_s^* \delta_0> ds 
\\
&\qquad
+ {\displaystyle\lim_{k\to\infty}}
 \frac 1{T_k} [\int_{T_k}^{t+T_k} <\phi,P_s^* \delta_0> ds
 -\int_1^{t+1} <\phi,P_s^* \delta_0> ds]
\\
&= <\phi,\mu> .
\end{array}
$$

As far as we are concerned with uniqueness of invariant measures, we
prove it along the lines of \cite{dpz2}. 
Let $P(t,y,\cdot)$ denote the transition probability: 
$P(t,y,\Gamma)= \mathbb P\{ u(t;y)\in \Gamma\}$.
By Khas'minskii theorem, if a Markovian semigroup $\{P_t\}_{t \ge 0}$ is
irreducible at time $t_1>0$
and strongly Feller at time $t_2>0$, then it is regular at time
$t_1+t_2$, that is the transition probabilities $P(t,y,\cdot)$
are equivalent for $t > t_1+t_2, y\in H$.
By Doob theorem, given an invariant measure $\mu$, 
if the Markovian semigroup $\{P_t\}_{t \ge 0}$ is regular for some
$t_0>0$,
then $\mu$ is strongly mixing and 
$$
 \lim_{t \to \infty} P(t,y,\Gamma)
= \mu(\Gamma)
$$
for arbitrary $y \in H$ and Borelian subset $\Gamma$  of $H$.
\\
Moreover, $\mu$ is the unique invariant measure, is ergodic 
 and is equivalent to any transition probability
measure $P(t,y,\cdot)$
for $y \in H$ and $t\ge t_0$.

Since the irreducibility and strongly Feller property are interesting
in themselves for a stochastic equation, we prove them separately in
Section \ref{irred} and \ref{strongly}, respectively.
First, we prove them in the space $H$; 
then, by regularity results in Section \ref{secGIR} we prove
them in more regular spaces.

\section{EXISTENCE OF INVARIANT MEASURES} \label{exinv}
In this section, we work with an operator $G$ satisfying
\eqref{ipoalfa} for some $\alpha>0$. 
\\
Keeping in mind that the space $D(A^{\tilde\alpha})$ is compactly embedded in
$H$ for any $\tilde\alpha>0$, our aim is to show that there exists a
parameter $\tilde\alpha
>0$ for which the following holds true:
\begin{equation}\label{tight}
 \forall \eps >0 \quad \exists R>0:\;
 \frac 1T \int_1^T\mathbb P \{ |A^{\tilde\alpha} u(t;0)|>R\} dt<\eps 
\mbox{ for all } T>1 .
\end{equation}
Notice that we consider initial data equal to zero and therefore the
regularity of the solution $u$ depends only on $G$. 

First of all, the linear equation enjoys  property \eqref{tight}.
Indeed, $z_a(t)=\int_0^t e^{-(\nu A^2-A+a)(t-s)}Gdw(s)$ and 
from \eqref{zinalfa}
we know that
$\displaystyle\sup_{0\le t<\infty} \mathbb E|A^\alpha z_a(t)|^2 $
is finite and tends to 0 as $a \to \infty$.
Moreover, by Chebyshev inequality
\begin{equation}\label{cheb}
\begin{split}
 \frac 1T \int_0^T \mathbb P\{|A^\alpha z_a(t)|>R\}dt 
&\le
 \frac 1T \int_0^T\frac {\mathbb E (|A^\alpha z_a(t)|^2) }{R^2}dt
\\
&\le
\frac 1{R^2} \sup_{0\le t<\infty}\mathbb E (|A^\alpha z_a(t)|^2) 
\end{split}
\end{equation}
and we conclude that $ \frac 1T \int_0^T \mathbb P\{|A^\alpha z_a(t)|>R\}dt$ 
can be made as small as we want by a
suitable choice of $R$.

Now, we prove
tightness for $u$, looking the equation satisfied by $u$
as a perturbation of the linear equation \eqref{eqOU}.
We follow \cite{dpg} (see also \cite{dpz2}).
As a first step, let us prove
\begin{proposition}\label{primoTIGHT}
Assume \eqref{ipotG}. Then 
\begin{equation}\label{tightH}
 \forall \eps >0 \quad \exists R>0:\;
 \frac 1T \int_0^T\mathbb P \{ |u(t;0)|>R\} dt<\eps \;\;
\mbox{ for all } T>0 .
\end{equation}
\end{proposition}
\proof
By \eqref{cheb}, this proposition
holds true if we prove \eqref{tightH} 
for $v_a=u-z_a$.

First, a priori estimates are required and we borrow from  the deterministic
case the suitable bounds.
The unknown $v_a$ satisfies \eqref{eq-v}
with $v_a(0)=0$.
We have proved that there exists a unique process solution such that
$v_a \in C([0,T];H)\cap L^2(0,T;D(A))$ $\mathbb P$-a.s. 
and
\begin{equation}\label{stima-v-a}
    \frac d{dt} |v_a|^2 +\nu |Av_a|^2 \le
   C_1 (1+|z_a|^2)|v_a|^2 + a^2 |z_a|^2+C_1|z_a|^4.
\end{equation}

Applying Gronwall lemma, we do not find useful estimates 
 to prove \eqref{tightH}. Hence, we proceed as in the deterministic
case.
We follow \cite{ce} (a similar result is in \cite{good}) 
and introduce an auxiliary function $h_s$, 
where $h_s(x)=h(x+s)$ for a suitable $h \in D(A)$ and
$s=s(t)$; other properties of $h$ will be presented below.
We work with $v_a-h_s\equiv v_a(t,x)-h(x+s(t))$; we have
\[
\begin{split}
\frac 12 \frac{d\;}{dt} |v_a-h_s|^2
&=
\langle \frac{d\;}{dt} v_a- \frac{d\;}{dt} h_s, v_a-h_s\rangle
\\
&= 
-\nu \langle A^2v_a,v_a\rangle
+\langle Av_a,v_a\rangle -\langle B(v_a),v_a-h_s\rangle
+\nu \langle Av_a,Ah_s\rangle
\\
&\quad
-\langle Av_a,h_s\rangle
+ a \langle z_a, v_a-h_s\rangle
-\langle B(v_a,z_a)+B(z_a,v_a),v_a-h_s\rangle
\\
&\quad
-\langle B(z_a),v_a-h_s\rangle
 -\langle h_s^\prime \dot s,v_a-h_s\rangle
\\
& \le
-\nu \langle A^2v_a,v_a\rangle
+\langle Av_a,v_a\rangle +\langle B(v_a),h_s\rangle
+\nu \langle Av_a,Ah_s\rangle
\\
&\quad
-\langle Av_a,h_s\rangle
+ a \langle z_a, v_a-h_s\rangle
-\langle B(v_a,z_a)+B(z_a,v_a),v_a-h_s\rangle
\\
&\quad
-\langle B(z_a),v_a-h_s\rangle ,
\end{split}
\]
where we used (see \cite{ce}) that $\langle h_s^\prime \dot s, v_a \rangle
=c \dot s^2 \ge 0$ and 
$\langle h_s^\prime,  h_s \rangle=0$.
\\
Generalizing  the proof of \cite{ce} (indeed, they consider $\nu =1$), 
we have that for any $\nu>0$ there exists a positive constant
$\kappa$, depending on $L$ and $\nu$,
such that
\[
-\nu \langle A^2v_a,v_a\rangle
+\langle Av_a,v_a\rangle +\langle B(v_a),h_s\rangle
\le
-\kappa \langle A^2v_a,v_a\rangle - \kappa |v_a|^2  .
\]
Therefore, using \eqref{scambi}, \eqref{buzu} and
Young inequality
\[
\begin{split}
\frac 12 \frac{d\;}{dt} |v_a-h_s|^2
+& \kappa \langle A^2 v_a,v_a \rangle 
+\kappa |v_a|^2
\\
&\le
\nu \langle Av_a,Ah_s\rangle
-\langle Av_a,h_s\rangle+a \langle z_a, v_a-h_s\rangle
\\
&\qquad
+2b(z_a,v_a,v_a-h_s)+b(v_a,v_a-h_s,z_a)
+\frac 12 b(z_a,v_a-h_s,z_a)
\\
&\le
\nu  |Av_a| |Ah_s| + |Av_a||h_s|+a|z_a||v_a-h_s|
\\
&\qquad
+c|z_a||Av_a||v_a-h_s|
+c|v_a||A(v_a-h_s)||z_a|
+c|z_a|^2|A(v_a-h_s)|
\\
&\le
\frac{\kappa} 3 |Av_a|^2 +\frac {\kappa}3  |A(v_a-h_s)|^2
+ c_\nu |z_a|^2|v_a-h_s|^2 +c_\nu |z_a|^2|v_a|^2
\\
&\qquad
+a|z_a||v_a-h_s|+c_\nu |z_a|^4
+c_\nu |Ah_s|^2 +c_\nu |h_s|^2 
\end{split}
\]
By means  of the triangle inequality and of embedding $D(A)\subset H$,
we get
that there exist positive constants $\kappa$ and $\Bar c$
(depending on $L$ and $\nu$) such that 
\begin{equation}\label{v-fi}
\frac{d\;}{dt} |v_a-h_s|^2 +\kappa |v_a-h_s|^2
\le
\Bar c |z_a|^2|v_a-h_s|^2
+\Bar c (1+|Ah|^4+|z_a|^4)
\end{equation}
since $|h_{s}|=|h|$.
\\
From \cite{dpg}, we define $\sigma(t)=\log(|v_a(t)-h_{s(t)}|^2\vee
R)$ 
for some
$R>1$ to be chosen
later. Denoting by $\chi_\Gamma$ the indicator function of the set $\Gamma$,
we have
\begin{equation}\label{intlog}
\int_0^T   \chi_{\{|v_a(t)-h_{s(t)}|^2\ge R\}}
 \frac {\frac{d\;}{dt}|v_a(t)-h_{s(t)}|^2 }{|v_a(t)-h_{s(t)}|^2}
 dt
\\=
\sigma(T)-\sigma(0) \ge 0.
\end{equation}
Multiplying \eqref{v-fi} by 
$\chi_{\{|v_a(t)-h|^2\ge R\}} \frac 1{|v_a(t)-h|^2} $, 
we obtain
\begin{multline*}
\chi_{\{|v_a(t)-h_{s(t)}|^2\ge R\}}
\frac { \frac{d\;}{dt}|v_a(t)-h_{s(t)}|^2 }{|v_a(t)-h_{s(t)}|^2}
+\kappa \chi_{\{|v_a(t)-h_{s(t)}|^2\ge R\}}
\\
\le
\Bar c |z_a(t)|^2+\frac {\Bar c}R (1+|Ah|^4+|z_a(t)|^4) .
\end{multline*}
Integrating in time, bearing in mind \eqref{intlog}  and taking expectation
we get
\begin{multline*}
 \frac {\kappa}T \int_0^T\mathbb P \{ |v_a(t)-h_{s(t)}|\ge R\} dt
\\
\le
 \frac {\Bar c}R(1+|Ah|^4) 
+ \frac {\Bar c}R \frac 1T \int_0^T\mathbb E |z_a(t)|^4 dt
+ \Bar c \sup_{0\le t< \infty}\mathbb E|z_a(t)|^2 .
\end{multline*}
Taking $a$ sufficiently large the last term can be as small as we
want;
moreover, for this fixed $a$
the two other terms in the r.h.s. tends to zero as $R$ tends
to $\infty$. On the other hand, by triangle inequality we find that
the same estimate holds for 
the quantity
$ \frac 1T \int_0^T\mathbb P \{ |v_a(t)|\ge R\} dt$.
Thus, the proposition is proved.  \hfill
$\Box$

\medskip

Now we prove \eqref{tight}.

\begin{proposition}
Let \eqref{ipoalfa} be satisfied for some $\alpha>0$. 
Then there exists $\tilde\alpha \in (0,1)$ such that
\begin{equation}\label{u-alfa}
 \forall \eps >0 \quad \exists R>0:\;
 \frac 1T \int_1^T\mathbb P \{ |A^{\tilde\alpha} u(t;0)|>R\} dt<\eps 
\mbox{ for all } T\ge 1. 
\end{equation}
Moreover, there exists at least one invariant measure for equation \eqref{sks}.
\end{proposition}
\proof
Since $A^{\alpha}
u= A^{\alpha} v_a +A^{\alpha} z_a$, 
according to \eqref{cheb} we seek the
result for $A^{\alpha} v_a$.
To deal with
 $A^\alpha v_a$, we exploit the regularizing effect 
of the semigroup
$e^{-\nu A^2 t}$.
Let us write equation \eqref{eq-v}, on the time interval $[t,t+1]$,  
in the integral form:
\[
\begin{split}
v_a(t+1)=e^{-\nu A^2}&v_a(t)+\int_t^{t+1} \!\!\!e^{-\nu A^2(t+1-s)}
    Av_a(s)\ ds 
\\
&
-\int_t^{t+1} \!\!\!e^{-\nu A^2(t+1-s)}B(u(s))\ ds
+a
\int_t^{t+1} \!\!\!e^{-\nu A^2(t+1-s)}z_a(s)\ ds .
\end{split}
\]
Then
\begin{multline}\label{Aalfa}
A^\alpha v_a(t+1)=A^\alpha e^{-\nu A^2}v_a(t)
+\int_t^{t+1} \!\!\!A^{\alpha +1} e^{-\nu A^2(t+1-s)}
    v_a(s)\ ds 
\\
-\int_t^{t+1} \!\!\!A^{\alpha +1} e^{-\nu A^2(t+1-s)}A^{-1}B(u(s))\ ds
+a
\int_t^{t+1} \!\!\!e^{-\nu A^2(t+1-s)}A^\alpha z_a(s)\ ds .
\end{multline}
By means of \eqref{regeA}, we estimate each term in the r.h.s.
\begin{align*}
& |A^\alpha e^{-\nu A^2}v_a(t)|\le K_1 |v_a(t)| ,
\\
& |\int_t^{t+1} \!\!\!A^{\alpha +1} e^{-\nu A^2(t+1-s)}
    v_a(s)\ ds|
\le 
\; \tilde K_2 
 \sup_{0\le r \le 1}|v_a(t+r)| ,
\intertext{where $\tilde K_2:=\displaystyle\int_t^{t+1}
     \frac{M_{\frac{\alpha+1}2}}{(t+1-s)^{\frac{\alpha+1}2}} ds$
is finite if $\alpha<1$. If $\alpha\ge 1$, we choose $\tilde
\alpha $ such that $0<\tilde \alpha <1 \le \alpha$ and obtain this
result for  $\tilde\alpha$ instead of $\alpha$.
}
&|\int_t^{t+1} \!\!\!A^{\alpha +1} 
 e^{-\nu A^2(t+1-s)}A^{-1}B(u(s))\
 ds|
\\&\qquad
\le
\int_t^{t+1} \frac{M_{\frac{\alpha+1}2}}
{(t+1-s)^{\frac{\alpha+1}2}} | B(u(s))|^2_{D(A^{-1})}ds
\\&\quad
\overset{\text{by }\eqref{tre}}{\le} 
\int_t^{t+1} 
 \frac {M_{\frac{\alpha+1}2}}{(t+1-s)^{\frac{\alpha+1}2}} c|u(s)|^2 ds
\\&\qquad
\le 
2c \sup_{0\le r \le 1} |v_a(t+r)|^2 \tilde K_2 
+ 2cM_{\frac{\alpha+1}2} \int_t^{t+1} \frac {|z_a(s)|^2}{(t+1-s)^{\frac{\alpha+1}2}} ds
\\&\qquad
\le
2 c\sup_{0\le r \le 1} |v_a(t+r)|^2 \tilde K_2 
+2cM_{\frac{\alpha+1}2} \tilde K_3 \big(\int_t^{t+1}|z_a(s)|^p
ds\big)^{2/p} ,
\end{align*}
where $\tilde K_3:=\Big(\int_t^{t+1}
     \frac{ds}{(t+1-s)^{\frac{\alpha+1}2 \frac p{p-2}}} \Big)
     ^{\frac{p-2}{p}}$ is finite
for suitable $p>2$ depending on $\alpha\in (0,1)$. Again,
 if $\alpha\ge 1$, we obtain the result for $\tilde
\alpha $ such that $0<\tilde \alpha <1 \le \alpha$.
The last bound is
\[
\begin{split}
|\int_t^{t+1} \!\!\!e^{-\nu A^2(t+1-s)}A^\alpha z_a(s)\ ds|
&\le
\int_t^{t+1} |A^\alpha z_a(s)|\ ds 
\\&
\le
\frac 12 \big(\int_t^{t+1}|A^\alpha z_a(s)|^p ds\big)^{2/p} +\frac 12 .
\end{split}
\]
Coming back to \eqref{Aalfa}, these bounds imply that there exists
$\tilde \alpha \in (0,1)$ such that
\[
\begin{split}
|A^{\tilde\alpha} v_a(t+1)|
\le 
K_1 |v_a(t)|
+ K_2\sup_{0\le r \le 1}|v_a(t+r)|^2
+K_3  \big(\int_t^{t+1}\!\!\!|A^\alpha z_a(s)|^p ds\big)^{2/p}
+K_4
\end{split}
\]
for suitable constants $K_1,K_2,K_3,K_4$.
Then
\[
\begin{split}
\mathbb P\{|A^{\tilde\alpha} v_a(t+1)|>R+K_4\}
\le
&\quad \mathbb P\{K_1 |v_a(t)|>\frac R3\}
\\
&+
\mathbb P\{K_2\sup_{0\le r \le 1}|v_a(t+r)|^2>\frac R3\}
\\
&+
\mathbb P\{ K_3\big(\int_t^{t+1}\!\!\!
         |A^\alpha z_a(s)|^p ds\big)^{2/p}>\frac R3\} .
\end{split}
\]
Bearing in mind  the proof of Proposition \ref{primoTIGHT}, we deal with the
first term in the r.h.s. in order to make it as small as we want by a
suitable choice of $R$. 
The same holds for the third term, using Chebyshev inequality and the
fact that, for $p/2$ integer,
$\mathbb E |A^\alpha z_a(s)|^p=c_p
[\mathbb E |A^\alpha z_a(s)|^2]^{p/2}$ by  Gaussianity. 
We are left with the second term to analyze.
We apply Gronwall lemma  to inequality
\eqref{stima-v-a} and  get
\[
\begin{split}
 \sup_{0\le r \le 1}|v_a(t+r)|^2 
& \le 
 |v_a(t)|^2 e^{\int_t^{t+1}
                      C_1(1+|z_a(s)|^2) ds}
\\
&\qquad 
+\int_t^{t+1} e^{\int_s^{t+1}  C_1(1+|z_a(r)|^2) dr}
 (a^2|z_a(s)|^2+C_1 |z_a(s)|^4)ds .
\end{split}
\]
Considering the probabilities, we have
\[ 
\begin{split}
\mathbb P\{\sup_{0\le r \le 1}&|v_a(t+r)|^2>R_2\} 
 \le
\mathbb P\{|v_a(t)|^2 e^{\int_t^{t+1}
                      C_1(1+|z_a(s)|^2) ds}>\frac {R_2}2\} 
\\
&\quad
+
\mathbb P\{\int_t^{t+1} e^{\int_s^{t+1}  C_1(1+|z_a(r)|^2) dr}
 (a^2|z_a(s)|^2+C_1|z_a(s)|^4)ds>\frac {R_2}2\} 
\\
&\le
\mathbb P\{|v_a(t)|^2>\sqrt{\frac {R_2}2}\} 
+
2 \mathbb P\{e^{\int_t^{t+1}
                      C_1(1+|z_a(s)|^2) ds}>\sqrt{\frac {R_2}2}\} 
\\
&\qquad
+
\mathbb P\{a^2|z_a(s)|^2+C_1|z_a(s)|^4>\sqrt{\frac {R_2}2}\}  .
\end{split}
\]
Once more, we deal with the first term in the last r.h.s according to
Proposition \ref{primoTIGHT} and with the two other terms by means of
Chebyshev inequality.
\\
This concludes the proof of \eqref{u-alfa}.
The existence of invariant measures follows from Krylov--Bogoliubov
method.  \hfill
$\Box$

\section{UNIQUENESS OF INVARIANT MEASURES}\label{unica}
We first prove irreducibility and strongly Feller property
in the space $H$. Then in Section \ref{secGIR}, we extend these
results
to more regular spaces $D(A^\alpha)$, $\alpha>0$,
 by means of Girsanov theorem.

\subsection{Irreducibility}\label{irred}
We say that 
the Markovian semigroup $\{P_t\}_{t\ge 0}$ is $H$-irreducible at time $t>0$
if, for arbitrary non empty open set $\Gamma\subseteq H$ and all $y \in H$,
$$
 P_t\chi_\Gamma(y)>0.
$$
This is equivalent to 
$$
 \mathbb P\{|u(t;y)-\tilde u|<R\}>0 \qquad
 \forall t>0 \quad\forall  y,\tilde u \in H\quad \forall R>0 .
$$
Given any $t>0$, we prove it 
following an idea from \cite{fm}: we show pathwise
that there are suitable $\bar R>0$ and $\bar z_a \in
C_0([0,T];H)$  
(the subset of $C([0,T];H)$-functions which vanish at $t=0$)
such that
\begin{equation}\label{differ}
 \{|u(t;y)-\tilde u|<R\} \supseteq \{|z_a-\bar z_a|_{C_0([0,t];H)}<\bar R\}
\end{equation}
and that the probability of the  last set is strictly  positive.

By \cite{dpz2} (see also \cite{masl}), 
if $G$ is a linear bounded operator in $H$ with range dense in $H$,
then 
the law of the process $z_a(\cdot;0)$ is full in $C_0([0,T];H)
$, i.e.
\[
 \mathbb P\{z_a\in \Lambda\}>0
\qquad \mbox{ for any non empty open set } \Lambda \subset
  C_0([0,T];H) .
\]
Notice that this covers the case of $G=A^\gamma$ 
for $\gamma\le 0$. But Lemma 2.6 in
\cite{masl} works also for $0<\gamma<\frac 34$; the hypothesis is
checked with calculations similar to
\eqref{zinalfa}. 

Hence, what remains to prove is \eqref{differ}.
First, we define a function $\bar u$ linking $y$ to $\tilde u$ in time $t$ as
$$
 \bar u(s)=y+\frac st [\tilde u -y ]\qquad \mbox{ for } 0\le s \le t.
$$
We have that $\bar u \in C([0,t];H)$.
Therefore
\begin{equation}\label{diff-u}
\{|u(t;y)-\tilde u|<R\} 
 \supseteq 
 \{|u(\cdot;y)-\bar u(\cdot;y)|_{C([0,t];H)}<R\} .
\end{equation}
We work now with the $C([0,t];H)$-norms and prove that
given $\bar u \in C([0,t];H)$ there exist $\bar v_a,
\bar z_a \in C([0,t];H)$ 
such that $\bar u = \bar v_a +\bar z_a$, $\bar v_a$ satisfies
equation \eqref{eq-v} with $z_a=\bar z_a$
and $\bar z_a(0)=0$; moreover 
\begin{equation}\label{diff-v}
 |v_a(\cdot;y)-\bar v_a(\cdot;y)|_{C([0,t];H)} \le 
 L |z_a -\bar z_a|_{C_0([0,t];H)}.
\end{equation}
Indeed, equation \eqref{eq-v}
satisfied by $v_a$ can be written as 
$$
 \left\{
\begin{array}{l}
 \frac{d}{dt} \bar v_a(t) + [\nu A^2 -A +a]\bar v_a(t)
   = a u(t) - B(\bar u(t))
 \\
 \bar v_a(0)=y
 \end{array}\right.
$$
By means of the estimates used in the previous sections, 
it is easy to check\footnote{We should work first with the Galerkin
  approximation $\bar v_a^n$ and then pass to limit as $n\to\infty$. But
we show the basic steps for $\bar v_a$.}
 that, given $\bar u \in C([0,t];H)$, the  r.h.s. 
$ a\bar u - B(\bar u)\in C([0,t];D(A^{-1}))$. 
Then, by classical results on
linear parabolic equations, we know that there exists
a unique solution $\bar v_a\in C([0,t];H)\cap L^2(0,t;D(A)) $.
\\
Therefore the function $\bar z_a := \bar u - \bar v_a$ is well defined and
belongs to $C_0([0,t];H)$.

The difference $V_a:=v_a-\bar v_a$ satisfies the equation
$$
\left\{
\begin{array}{ll}
 \frac{d\;}{dt} V_a + \nu A^2 V_a &-AV_a +B(V_a,v_a)+
  B(\bar v_a,V_a)+B(V_a,z_a)+B(\bar v_a,Z_a)
 \\
&+B(Z_a,v_a)+B(\bar z_a, V_a)=AZ_a -B(Z_a,z_a)-B(\bar z_a, Z_a)
 \\
 V_a(0)=0
 \end{array}\right.
$$
(where $Z_a=z_a-\bar z_a$).

Taking the scalar product in $H$ of the first equation with $V_a$ and
using
estimates on the trilinear form as usual, we get
\[
\begin{split}
 \frac 12 \frac {d}{dt} |V_a|^2+\nu |AV_a|^2
&\le
   c(|v_a|+|\bar v_a|+|z_a|+|\bar z_a|) |V_a||AV_a|
   +c|A \bar v_a| |V_a| |Z_a| 
\\
&
  \qquad +  c(1+|\bar v_a|+|z_a|+|\bar z_a|)|Z_a||AV_a|
\\
& 
  \le \frac \nu 2 |AV_a|^2 + \psi_1 |V_a|^2 + \psi_2 |Z_a|^2
\end{split}
\]
where $\psi_1=c_\nu(1+|v_a|^2+|\bar v_a|^2+|z_a|^2+|\bar z_a|^2) \in
C([0,t];H)$
and $\psi_2=c_\nu(1+|A \bar v_a|^2+|z_a|^2+|\bar z_a|^2)
\in L^1(0,t;H)$. 
 
By Gronwall lemma, we conclude that
$$
 \sup_{0\le s\le t} |V_a(s)|^2 \le 
 \int_0^t e^{\int_s^t \psi_1(r)dr}\psi_2(s)\, |Z_a(s)|^2 ds .
$$
Hence
$$
 \sup_{0\le s\le t} |V_a(s)|^2 \le L \sup_{0\le s\le t} |Z_a(s)|^2 
$$
for some constant $L$ depending on $t$, on the $H$- norm of $y$ and
the $C([0,t];H)$-norm of $z_a$, $\bar z_a$.
This is \eqref{diff-v};
it implies that
$$
 |u(\cdot;y)-\bar u(\cdot;y)|_{C([0,t];H)} \le 
 (1+L) |z_a -\bar z_a|_{C_0([0,t];H)}
$$ 
and by \eqref{diff-u}, 
the relationship  \eqref{differ} follows with $ R=(1+L)\bar R$.

Summing up, we have proved the following result.
\begin{proposition}\label{strong}
Let $G$ be either a linear bounded operator in $H$
satisfying \eqref{ipotG} and  with dense range in $H$ or $G=A^\gamma$
with $0<\gamma<\frac 34$.
Then the Markovian semigroup $\{P_t\}_{t \ge 0}$ is $H$-irreducible
for any $t>0$.
\end{proposition}

\begin{remark}
If $G$ is diagonal, then all the components have to be non zero,
i.e. $Gw(t)=\sum_{j=1}^\infty g_j \beta_j(t)e_j $ with all $g_j \neq
0$
and $\sum_j g_j^2 \lambda_j^{-2}<\infty$.
\\
Of course, the result holds also for $G=LA^\gamma$
with $\gamma<\frac 34$ and $L$ an isomorphism in $H$. 
\end{remark}

\subsection{Strongly Feller}\label{strongly}

The Markovian semigroup $\{P_t\}_{t\ge 0}$
is strongly Feller in $H$ at time  $t>0$ if
$$
 P_t : B_b(H) \to C_b(H).
$$ 
This means that, given $\phi \in B_b(H)$
$$
 \lim_{y_1 \to y_2} P_t\phi(y_1)=P_t\phi(y_2).
$$

If the operator $G$ fulfils  the assumptions of Proposition
\ref{strong} below, 
i.e. if $G$ is regular enough and invertible, it is well-known 
that the Markovian semigroup for the linear equation
\eqref{eqOU} is strongly Feller for any $t>0$ (see e.g. \cite{dpz2,masl}).
Considering equation \eqref{sks} 
as a non linear version of \eqref{eqOU},
the Markovian semigroup $P_t$ may have this
regularizing effect for
$t>0$. 
There is a formula for its derivative, providing
even more regularity under some assumptions on the non linear part 
(see \cite{dpz2}).
To use it, we need to 
introduce a modified Galerkin equation; then in the limit we recover
our equation and the strongly Feller property will be proved.
Our technique is similar to that of  \cite{fm,fe}. 

For $R\ge 1$, let $\theta_R$ be a $C^1$-function with bounded derivative
such that $\theta_R=1$ on $[-R,R]$ and $\theta_R=0$
outside $[-R-1,R+1]$; take $|\theta_R|$ and $|\theta_R^\prime|$ bounded by 1. 
We consider the modified Galerkin system
\begin{equation}\label{thetaGal}
\begin{cases}
 du^{n,R}+\big[\nu A^2 u^{n,R}
   -Au^{n,R}+\theta_R(|u^{n,R}(t)|^2)B_n(u^{n,R})\big]\,dt=\Pi_n Gdw
\\
u^{n,R}(0)=\Pi_n y
\end{cases}
\end{equation}
The cut-off function $\theta_R$ introduces minor changes with respect
to the Galerkin equation for \eqref{sks}; it is straightforward 
to obtain, as in Section 3, that there exists a unique process $u^{n,R}
\in C([0,T];H)$ solving \eqref{thetaGal}.

First we show that the Markovian semigroup $\{P^{n,R}_t\}_{t\ge 0}$ 
of  \eqref{thetaGal}
is Lipschitz Feller for $t>0$.
Let $[D P^{n,R}_t\phi(y)]\cdot h$ be the derivative,  in direction
$h$,
 of the mapping $y \mapsto P^{n,R}_t\phi(y)$. 
Then, Bismut-Elworthy-Li formula 
represents it by means of an expression depending on 
$[Du^{n,R}(t;y)]\cdot h$, the derivative, in the direction $h$, 
of the mapping $y \mapsto u^{n,R}(t;y)$.
Indeed (see \cite{dpz2} and references therein)
\begin{equation}\label{bh}
[D\: P^{n,R}_{t} \phi (y)] \cdot h 
 = 
{1 \over t} \mathbb E \left[ \phi(u^{n,R}(t;y)) \int_{0}^{t}
        \langle (\Pi_{n}GG^{*}\Pi_{n})^{-{1 \over 2}} 
         [D\: u^{n,R}(s;y)] \cdot h,
        \Pi_n dw(s)\rangle \right]
\end{equation}
for all $h \in H_{n}$.
Therefore
$$
 \left|D\: P^{n,R}_{t} \phi (y) \cdot h \right| 
\le
 {1 \over t} \| \phi \|_0   
        \left[ \mathbb E \int_0^t 
         \left|(\Pi_{n}GG^{*}\Pi_{n})^{-{1 \over 2}} 
        D\: u^{n,R}(s;y) \cdot h\right|^{2} ds  \right]^{1 \over 2}.
$$
Here $\| \cdot \|_{0}$ denotes the supremum norm in $C_{b}$ (or $B_{b}$).
\\
Assuming that $G$ is a linear bounded operator in $H$ with
$R(G) \supseteq D(A)$,
we have that
$$
\int_0^t 
         \left|(\Pi_{n}GG^{*}\Pi_{n})^{-{1 \over 2}} 
        D\: u^{n,R}(s;y) \cdot h\right|^{2} ds 
\le C
\int_0^t 
         \left| A 
        D\: u^{n,R}(s;y) \cdot h\right|^2 ds 
$$
(see details in \cite{fm}). This holds also 
if $G=LA^\gamma$ with $0<\gamma<\frac 34$ as in Remark \ref{convergenza}.\\
Estimate on the latter quantity are easily obtained; by
\eqref{thetaGal}, the equation for $U^{n,R}:=D\: u_{n}^{R}(\cdot;y) \cdot h$
is
$$
\begin{cases}
 \frac{d\;}{dt} U^{n,R} +\nu A^2 U^{n,R}-A U^{n,R} 
 +2
 \theta_R^\prime(|u^{n,R}|^2) \langle u^{n,R},U^{n,R}\rangle B_n(u^{n,R})
\\
\qquad\qquad\qquad\qquad +
 \theta_R(|u^{n,R}|^2) [ B_n(u^{n,R},U^{n,R})+B_n(U^{n,R},u^{n,R})]=0 ,
\\
U^{n,R}(0)=\Pi_n h.
\end{cases}
$$
We multiply scalarly in $H$ the first equation by $U^{n,R}$; by means
of relationships for the trilinear forms already used before, we obtain
$$
 \frac 12 \frac{d\;}{dt} |U^{n,R}|^2+\nu |AU^{n,R}|^2 
\le
 \frac \nu 2 |AU^{n,R}|^2
 +c_\nu \big[1+(\theta_R^\prime |u^{n,R}|^3)^2+
  (\theta_R |u^{n,R}|)^2\big] |U^{n,R}|^2 .
$$
Then
$$
 \frac{d\;}{dt} |U^{n,R}|^2+\nu |AU^{n,R}|^2 
\le 
  \tilde C |U^{n,R}|^2
$$
for some constant $\tilde C$ dependent on $\nu$ but not on $R$ or $n$.

By Gronwall lemma, 
$$
 \mathbb E |U^{n,R}(t)|^2 \le |h|^2 e^{\tilde C t}
$$
and, integrating in time,
\begin{equation}\label{stimaAU^R}
 \mathbb E \int_0^t |AU^{n,R}(s)|^2 ds \le \frac 1\nu |h|^2 (1+e^{\tilde
 C t}) .
\end{equation}
Summing up,
we have shown that
$$
 \left|D\: P^{n,R}_{t} \phi (y) \cdot h \right| 
 \le
  {1 \over t} \| \phi \|_0  c_\nu (1+e^{\tilde C t})^{1/2} |h|
=: L_{R,t} \| \phi \|_0  |h| .
$$
Thus, the derivative $D\: P^{n,R}_{t} \phi (y)$ is uniformly bounded 
and by the  mean value Theorem
$$
\left|P^{n,R}_{t} \phi (y_1) - P^{n,R}_{t} \phi (y_2) \right| 
        \le 
\: \Big(\sup_{ \stackrel{k,h \in H_{n}}{ |h| \le 1}} 
        \left|D P^{n,R}_{t} \phi (k) \cdot h\right| \Big)\; |y_1-y_2|
\le 
L_{R,t}  \|\phi\|_0  \; |y_1-y_2|       ,
$$
that is  $P^{n,R}_{t}\phi$ is Lipschitz Feller for all $t>0$.

Passing to the limit, as $n \to \infty$, we obtain
as usual that, $\mathbb P$-a.s.,  
$u^{n,R}$ converges strongly to $u^R$ in $L^2(0,T;H)$, where 
$u^R$ solves the equation
$$
  du^R+\big[\nu A^2 u^R
   -Au^R+\theta_R(|u^R|^2)B(u^R)\big]\,dt=Gdw
$$
Passing to a subsequence, $u^{n,R}(t)$ converges to 
$u^R(t)$ in $H$, for a.e. $t$.
Then, for $\phi \in C_b(H)$,  a subsequence of  $P^{n,R}_t\phi(y)$
converges  towards
$P^{R}_t\phi(y)$, for a.e. $t$. But the trajectories of $u^R$ are
continuous in time with values in $H$ and therefore we conclude that
for any $t>0$, $R\ge 1$, $\phi \in C_b(H)$ and $y_1,y_2 \in H$
there exists a constant $L_{R,t}$ depending only on $R$ and $t$, such that
$$
\left|P^{R}_{t} \phi (y_1) - P^{R}_{t} \phi (y_2) \right| 
    \le 
L_{R,t}  \|\phi\|_0 \; |y_1-y_2|.  
$$
The same result holds for $\phi \in B_b(H)$ (see 
Lemma 7.1.5 in \cite{dpz2}).

The last step consists in letting $R \to \infty$.
Working as we did for $v_n$ obtaining  \eqref{sup-v}, we can
easily verify that
$$
 \sup_{|y|\le M} \; \sup_{0\le t \le T}|u(t;y)|< \infty 
\qquad \mathbb P-a.s.
$$
and similarly for $u^R(t;y)$.
Moreover, 
the processes $u$ and $u^R$ coincide until $u^R$ lies 
in the ball of radius $R$ in $H$. 
Therefore, given $t>0$ and $y \in H$, for $\mathbb
P$-a.e. $\omega$ there exists $R_\omega$ such that
$u^R(t;y)(\omega)=u(t;y)(\omega)$ for all $R\ge R_\omega$, uniformly
in $y$ in bounded sets of $H$.
So, given $t$ and $\phi \in B_b(H)$, $\mathbb P$-a.s.
$$
 \lim_{R\to \infty} \phi(u^R(t;y))=\phi(u(t;y))
$$
uniformly in $y$ in bounded sets of $H$.
Since $\phi$ is bounded, we get that the convergence holds when we
take expectation.
Hence, for any $t>0$ and $\phi\in B_b(H)$
$$
 \lim_{R\to \infty} P_t^R\phi(y)=P_t\phi(y)
$$
uniformly in $y$ in bounded sets of $H$.

Finally, given $t>0$, $\phi \in B_b(H)$ and $y_1,y_2\in H$ 
\[
\begin{split}
\lim_{y_1\to y_2} P_t \phi(y_1)
&=
\lim_{y_1\to y_2} \lim_{R\to \infty} P^{R}_t\phi(y_1)
=
\lim_{R\to \infty} \lim_{y_1\to y_2} P^{R}_t\phi(y_1)
=
\lim_{R\to \infty} P^{R}_t\phi(y_2)
\\
&=
P_t \phi(y_2).
\end{split}
\]
This proves the strongly Feller property for any $t>0$.
Therefore, we have proved the following result.
\begin{proposition}
Let $G$ be either a linear bounded operator in $H$, satisfying
\eqref{ipotG}
and such that  $R(G)\supseteq D(A)$
or $G=A^\gamma$
with $0<\gamma<\frac 34$.
Then, for any $y \in H$
 the Markovian semigroup $\{P_t\}_{t \ge 0}$ is 
strongly Feller
in $H$ for any $t>0$.
\end{proposition}

\begin{remark}\label{ossSF}
For instance, $G=A^\gamma$, with $-1 \le \gamma<\frac 34$, 
is an example of operator fulfilling the assumptions of
the previous proposition.
The case $G=LA^\gamma$, with $L$ an isomorphism in $H$, 
can be treated in the same way.
\end{remark}

\subsection{Regularity results}\label{secGIR}
In this section we work with the operator $G$ of the form $A^\gamma$,
considering the values  $\gamma <-1$ not included in the previous
section (see Remark \ref{ossSF}).
 Even if the problem interesting from the physical point
of view (that with $\gamma =\frac 12 $)
has been solved in the previous sections, we want to show that the
limitation $\gamma \ge -1$ can be removed so to prove existence and
uniqueness of the invariant measure when $G=A^\gamma$ for any
$\gamma<\frac 34$.
The case $LA^\gamma$, where $L$ is an isomorphism in $H$, can be treated in
the same way.

We prove the following result.
\begin{proposition}\label{prorego}
Let $G=A^\gamma$ with $\gamma <-1$. For any $\alpha$ such that
\begin{equation}\label{condiz}
 1 \le \alpha < \frac 34 -\gamma,
\end{equation}
given $y \in D(A^\alpha)$
there exists a unique solution $u$ to problem \eqref{sks} and 
$$
 u \in C([0,T];D(A^\alpha)) \; \mathbb P-a.s.
$$
\end{proposition}
\proof
Existence and uniqueness hold in the bigger space $C([0,T];H)$, as proved in
Theorem \ref{primoTEO}. 
Thus, we need to prove the regularity
$C([0,T];D(A^\alpha))$. From \eqref{ipoalfa}, which now reads 
$\sum_j \lambda_j ^{2(\gamma+\alpha-1)}<\infty$,  
we obtain  that if $\alpha <\frac 34 -\gamma$ the linear equation has
solution
$z_a$ with paths in $C([0,T];D(A^\alpha))$.
What remain to be proved is that equation \eqref{eq-v} has a solution $v_a 
$ with paths in  $C([0,T];D(A^\alpha))$.
We show a priori estimates as in the proof of Proposition \ref{vinH}.
We take the scalar product of this equation with $A^{2\alpha}v_a$,
use \eqref{sei} and Young inequality  to  obtain
\[
\begin{split}
 \frac 12 \frac d{dt} |A^\alpha v_a|^2 &+\nu |A^{1+\alpha}v_a|^2 
\\
&= \langle A^\alpha v_a, A^{1+\alpha}v_a\rangle
   -\langle A^{\alpha-1} B(v_a+z_a),   A^{1+\alpha}v_a \rangle
   +a \langle A^\alpha z_a, A^\alpha v_a\rangle  
\\
&\le 
 |A^\alpha v_a||A^{1+\alpha}v_a| +
  | A^{\alpha-\frac 12} (v_a+z_a)|^2 |A^{1+\alpha}v_a|
 + a  |A^\alpha z_a| | A^\alpha v_a|
\\
&\le
 \frac \nu 2 |A^{1+\alpha}v_a|^2
+ \frac {c_\nu}2
 [|A^\alpha v_a|^2+|A^\alpha z_a|^2+| A^{\alpha-\frac 12} (v_a+z_a)|^4].
\end{split}
\]
Hence
\begin{equation}\label{stime-reg}
\frac d{dt} |A^\alpha v_a|^2+\nu |A^{1+\alpha}v_a|^2 
\le
  c_\nu |A^\alpha v_a|^2
 +c_\nu[|A^\alpha z_a|^2+| A^{\alpha-\frac 12} (v_a+z_a)|^4].
\end{equation}
First, consider  $\alpha=1$; then 
$$
 \frac d{dt} |A v_a|^2+\nu |A^{2}v_a|^2
\le c_\nu |A v_a|^2+c_\nu[|A z_a|^2+| A^{\frac 12} (v_a+z_a)|^4].
$$
But
$| A^{\frac 12} (v_a+z_a)|^4=
|A^{\frac 12} u |^4$ and it belongs to $ L^1(0,T)$; 
indeed, by Corollary \ref{cor-alfa} (ii)
we know that $u\in C([0,T];H)\cap L^2(0,T;D(A))$ and by interpolation 
$u \in L^4(0,T;D(A^{\frac 12}))$.
Since $|A z_a|^2+| A^{\frac 12} (v_a+z_a)|^4$ is  integrable in time, 
by Gronwall lemma we get 
as usual that
$$
 \sup_{0\le t \le T} |A v_a(t)|<\infty \, , \qquad
 \int_0^T |A^{2}v_a(t)|^2 dt <\infty.
$$
We do not consider the Galerkin approximations but we should do, as in
Section \ref{soluts}, to conclude that
$$
 v_a \in C([0,T];D(A))\cap L^2(0,T;D(A^2)).
$$
All the results hold pathwise, as before, 
 and we do not specify this every time.

Coming back to \eqref{stime-reg}, if $1<\alpha \le \frac 32$, then
\[
\frac d{dt} |A^\alpha v_a|^2+\nu |A^{1+\alpha}v_a|^2 
\le c_\nu |A^\alpha v_a|^2+c_\nu[|A^\alpha z_a|^2+| A (v_a+z_a)|^4]
\]
From the case $\alpha=1$ we know that $u=v_a+z_a \in C([0,T];D(A))$
and therefore $|A^\alpha z_a|^2+| A (v_a+z_a)|^4$ belongs to
$L^1(0,T)$.
Once more, Gronwall lemma allows to conclude that
$$
  \sup_{0\le t \le T} |A^\alpha v_a(t)|<\infty \, , \qquad
 \int_0^T |A^{1+\alpha}v_a(t)|^2 dt <\infty
$$
for $1<\alpha \le \frac 32$
and therefore
$
 v_a \in C([0,T];D(A^\alpha))\cap L^2(0,T;D(A^{1+\alpha}))
$.

By induction, we prove the result for all $\alpha \ge 1$.
Assume that for some integer $m>2$, given 
$\frac m2< \alpha\le \frac
{m+1}2$ we have $v_a \in C([0,T];D(A^\alpha))\cap
L^2(0,T;D(A^{1+\alpha}))$. 
Then we want to show the regularity result  for  $\frac
{m+1}2 <\alpha \le \frac {m+2}2$.
First, for $\alpha \le \frac {m+2}2$ the term
$| A^{\alpha-\frac 12} (v_a+z_a)|^4$ in  \eqref{stime-reg}
is bounded by
$| A^{\frac {m+1}2} (v_a+z_a)|^4$. But,
by the induction assumption we know in particular that
$u=v_a+z_a \in C([0,T];D(A^{\frac {m+1}2}))$.
Hence, as before,  we have the suitable a priori estimate 
$$
\frac d{dt} |A^\alpha v_a|^2+\nu |A^{1+\alpha}v_a|^2 
\le c_\nu |A^\alpha v_a|^2+c_\nu[|A^\alpha z_a|^2
+| A^{\frac {m+1}2} (v_a+z_a)|^4]
$$
for $\frac {m+1}2 <\alpha \le \frac {m+2}2$, which
allows to conclude by Gronwall lemma that
$$
 \sup_{0\le t \le T} |A^\alpha v_a(t)|<\infty \, , \qquad
 \int_0^T |A^{1+\alpha}v_a(t)|^2 dt <\infty
$$
so that $v_a \in C([0,T];D(A^\alpha))\cap
L^2(0,T;D(A^{1+\alpha}))$ for $\frac {m+1}2 <\alpha \le \frac {m+2}2$.
 \hfill
$\Box$

\medskip

Now, with the same estimates we can prove
as in  Section \ref{strongly}
that the Markovian semigroup is strongly Feller in $D(A^{\alpha})$, 
 with $\alpha$ specified by \eqref{condiz}.
The only differences with respect to Section
\ref{strongly} are that the cut-off function in \eqref{thetaGal}
is $\theta_R(|A^\alpha
u|^2)$ instead of $\theta_R(|u|^2)$ and we require
$R(G) \supseteq D(A^{1+\alpha})$
instead of $R(G) \supseteq D(A)$. 
This condition on the range of $G=A^\gamma$ is
satisfied if 
$$
1+\alpha \ge -\gamma.
$$

Moreover, we notice that Proposition \ref{prorego}
implies that 
the Markovian semigroup associated to equation \eqref{sks} 
is irreducible 
in $D(A^\alpha)$. Indeed, in Proposition \ref{strong} we  proved, for any
$\gamma<\frac 34 $, 
irreducibility in $H$. 
On the other hand, by Proposition \ref{prorego} we know that
the Markov process $u$ lives in $D(A^\alpha)$ for $\alpha$ specified
by \eqref{condiz}. Hence, the
irreducibility in inherited from $H$, since 
an open non empty  subset of $H$, when restricted to
 $D(A^{\alpha})$, is an open non empty  subset  of $D(A^{\alpha})$.

Summing up all the conditions on $\alpha$, we have the following result
\begin{proposition}
Given $\gamma<-1$, the Markovian semigroup $\{P_t\}$
associated to equation \eqref{sks} is irreducible and strongly Feller
in $D(A^\alpha)$ for any $t>0$, if $\alpha$ satisfies
\begin{equation}\label{alfa-gam}
 \begin{cases}
   1 \le  \alpha <\frac 34 -\gamma & \text{ when } -2 \le \gamma<-1\\
   -1-\gamma \le  \alpha <\frac 34 -\gamma & \text{ when } \gamma<-2
 \end{cases}
\end{equation}
\end{proposition}

\section{CONCLUSIONS}\label{conclu}
We collect all the previous results which imply existence 
of an invariant measure and
its uniqueness, as explained in Section \ref{invmeas}.

\begin{theorem}\label{concludo} 
Let $G=LA^\gamma$ for $\gamma <\frac 34$ and $L$ an isomorphism in $H$.
Set $E=H$ if $\gamma \ge -1$ and 
$E=D(A^\alpha)$ if $\gamma<-1$, where $\alpha$ fulfills 
\eqref{alfa-gam}.
Then equation \eqref{sks} has  a unique invariant measure $\mu$, 
concentrated on the space $E$.
All the transition probability
measures $P(t,y,\cdot)$
for $y \in E$ and $t>0$ are equivalent to $\mu$
and this measure $\mu$
is ergodic:
$$
 \lim_{T\to \infty}\frac 1T \int_0^T \phi(u(t;y)) dt 
  =\int \phi d\mu
$$
$\mathbb P$-a.s. and 
for all $\phi \in L^1(E;\mu)$, $y \in E$.
\end{theorem}

\end{document}